\documentclass[12pt]{article}
\usepackage[amsmath]{e-jc}
\usepackage{graphicx}
\theoremstyle{plain}
\newtheorem{axiom}[theorem]{Axiom}
\dateline{Dec 22, 2022}{TBD}{TBD}
\Copyright{The author. Released under the CC BY-ND license (International 4.0).}
\title{A Result on the Small Quasi-Kernel Conjecture}
\author{Allan van Hulst \\
\small VU University Amsterdam \\
\small The Netherlands \\
\small\tt a.c.van.hulst@vu.nl}
\begin{document}
\maketitle
\begin{abstract}
Any directed graph $D=(V(D),A(D))$ in this work is assumed to be finite and without 
self-loops. A source in a directed graph is a vertex having at least one ingoing 
arc. A quasi-kernel $Q\subseteq V(D)$ is an independent set in $D$ such that every 
vertex in $V(D)$ can be reached in at most two steps from a vertex in $Q$. It is 
an open problem whether every source-free directed graph has a quasi-kernel 
of size at most $|V(D)|/2$, a problem known as the small 
quasi-kernel conjecture (SQKC). The aim of this paper is to prove the SQKC 
under the assumption of a structural property of directed graphs. This relates 
the SQKC to the existence of a vertex $u\in V(D)$ and a bound on the number of 
new sources emerging when $u$ and its out-neighborhood are removed from $D$. The
results in this work are of technical nature and therefore additionally 
verified by means of the Coq proof-assistant.
\end{abstract}

\section{Introduction}
\label{sec:intro}
Certain types of independent sets in directed graphs can be used for a variety of
applications. For example, they can be used to model how a contagious disease spreads 
from a set of independent originators further into a population. Or, how a news message 
is communicated over a network if a set of independent sources starts spreading 
the message. 

Informally, a quasi-kernel is an independent set in a directed graph such that 
every vertex can be reached in at most two steps from this quasi-kernel. The main 
aim of this paper is to prove a conditional result about quasi-kernels. More 
specifically, to prove that every source-free directed graph has a quasi-kernel of size 
at most half of the number of vertices. This relies on the assumption of a structural 
property of directed graphs, which is stated as an axiom in this work.

The existence problem of such a `small' quasi-kernel (i.e. of size at most half of 
the number of vertices) has a long history, dating back to 1976, although the earliest 
published reference is quite recent \cite{fete}. It is originally attributed to 
P.L. Erd\H{o}s and L.A. Sz\'ekely \cite{lemon} and became known as the `small 
quasi-kernel conjecture'.

Only partial solutions are known at the moment. Disjoint quasi-kernels have been studied 
before but cannot be guaranteed to exist \cite{heard,gutin}. Presence 
of a kernel implies the existence of a quasi-kernel (\cite{arxiv}, further simplified in 
\cite{improved}), but many directed graphs do not have kernels. A slightly alternative 
approach was recently proposed, which suggests a bound in terms of the number of sources 
\cite{kostochka}.

The starting point of the technique applied in this paper is a well-known
constructive lemma which shows that every directed graph has a 
quasi-kernel \cite{chvatal}. Based on the recursive way the directed graph is broken 
down in this lemma, a technique is derived to count the number of vertices in the
quasi-kernel. If the number of new sources introduced by each of these steps
can be limited, this indeed results in a quasi-kernel of at most the desired 
bound.

The remainder of this paper is organized as follows. Section \ref{sec:def} lists
a number of definitions and in Section \ref{sec:breakdown} the way in which directed
graphs are broken down is formalized. Section \ref{sec:main} then details a proof
of the small quasi-kernel conjecture, under the axiom described previously. 

As a final introductory remark, the techniques proposed in this work are of 
somewhat technical nature. However, this is necessary as no alternative
formalism seems to be available for such a construction. All definitions and
proofs have been formalized by means of the Coq proof-assistant in order to
ensure absolute correctness. However, the paper is self-contained and one can
verify all results without perusing the accompanying Coq script. The Coq code
is also available through a GitHub repository:

\begin{center}
\href{https://github.com/allanvanhulst/sqkc/blob/main/qkernel.v}
     {\texttt{https://github.com/allanvanhulst/sqkc/blob/main/qkernel.v}}
\end{center}

\begin{figure}
\begin{center}
\includegraphics[scale=.75]{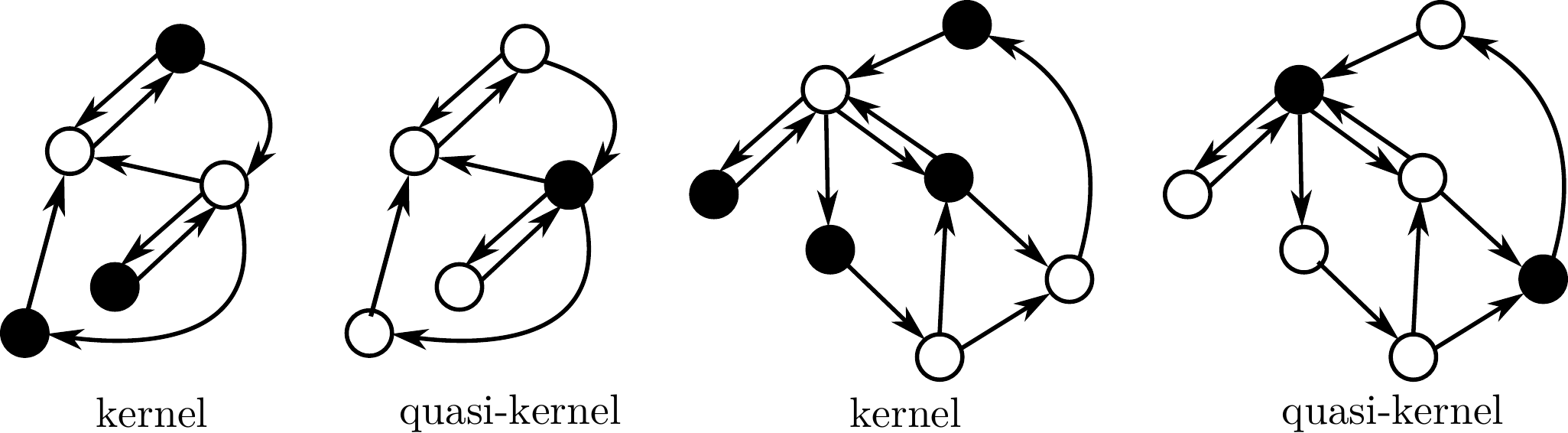}
\end{center}
\caption{Two examples of directed graphs with kernels and quasi-kernels, drawn as
  black vertices. In both graphs all possible kernels are larger than half of the 
  number of vertices.}
\label{fig:first}
\end{figure}

\section{Definitions}
\label{sec:def}
Any directed graph $D=(V(D),A(D))$ considered here is assumed to be finite and without 
self-loops. The notation $u\longrightarrow v$ is used to denote that $(u,v)\in A(D)$.
The open out-neighborhood (notation: $N^+(u)$) and closed out-neighborhood ($N^+[u]$) 
of a vertex $u\in V(D)$ are defined as:
\begin{center}
\begin{math}
N^+(u) = \{v\mid u\longrightarrow v\}\qquad
N^+[u] = N^+(u)\cup\{u\} 
\end{math}
\end{center}
The out-degree $d^+(u)$ for a vertex $u\in V(D)$ is then defined simply as 
$|N^+(u)|$. 

A vertex $u\in V(D)$ is said to be a source in $D$ if it does not have any ingoing
arc. Note that this classifies an isolated vertex as a source. Furthermore, a directed 
graph is said to be source-free if it does not have any sources (i.e. every vertex has 
at least one ingoing arc). Examples of source-free directed graphs are shown in Figure 
\ref{fig:first}.

If $D$ is a directed graph and $S\subseteq V(D)$ then the notation $D-S$ will be used 
to denote the subgraph induced by $V(D)-S$, formally:
\begin{center}
\begin{math}
D - S = (V(D)-S,\{(u,v)\in A(D)\mid u\not\in S,\,v\not\in S\}).
\end{math}
\end{center}
\noindent Informally, this is the directed subgraph where all arcs incident to a 
vertex in $S$ are removed from $D$.

Kernels in directed graphs were introduced by Richardson \cite{richardson}. 
A kernel $K\subseteq V(D)$ is an independent set in $D$ such that for every vertex 
$v\in V(D)$ it holds that either $v\in K$ or there exists a $u\in K$ such that 
$u\longrightarrow v$. Note that there are simple directed graphs without a kernel,
for instance a directed triangle.

A quasi-kernel is a weakening of the concept of kernel and defined as follows. A
quasi-kernel $Q\subseteq V(D)$ is an independent set in $D$ such that for all $w\in V(D)$ 
at least one of the following conditions is satisfied: (1) $w\in Q$ or, (2) there exists
a $v\in Q$ such that $v\longrightarrow w$ or, (3) there exist $u\in Q$ and $v\in V(D)$ such
that $u\longrightarrow v$ and $v\longrightarrow w$. Informally, a quasi-kernel $Q$ is 
therefore an independent set in $D$ such that every vertex in $V(D)$ can be reached in
at most two directed steps from $Q$. It is a well-known result that every directed
graph has at least one quasi-kernel (\cite{chvatal}, see also Lemma 
\ref{lem:exists}). Examples of kernels and quasi-kernels are shown in Figure 
\ref{fig:first}.

Clearly, every kernel is also a quasi-kernel but the examples in Figure \ref{fig:first} 
show these do not necessarily have to be under the desired bound of $|V(D)|/2$. However,
it was proved recently that the presence of a kernel is sufficient to derive the 
existence of a quasi-kernel which indeed satisfies this bound \cite{arxiv}.

For the counting method proposed in this paper, it will be necessary to relate
vertices in a quasi-kernel to out-neighbors of a specific type. This is a slight
adaptation of a similar concept in the theory of dominating sets in undirected graphs 
\cite{cockayne}.

\begin{definition}
\label{def:epon}
An \underline{external private out-neighbor} (\textsc{epon}) 
of a vertex $u\in S$, with regard to a subset $S\subseteq V(D)$, is another vertex 
$v\in V(D)$ such that $u\longrightarrow v$ and $v\not\in S$, with the additional 
property that $u$ is the only in-neighbor of $v$ in $S$. 
\end{definition}

\noindent For example, all directed graphs, except the first, in Figure \ref{fig:first} 
have \textsc{epon}s with regard to the subsets drawn as black vertices.

\begin{figure}
\begin{center}
\includegraphics[scale=.75]{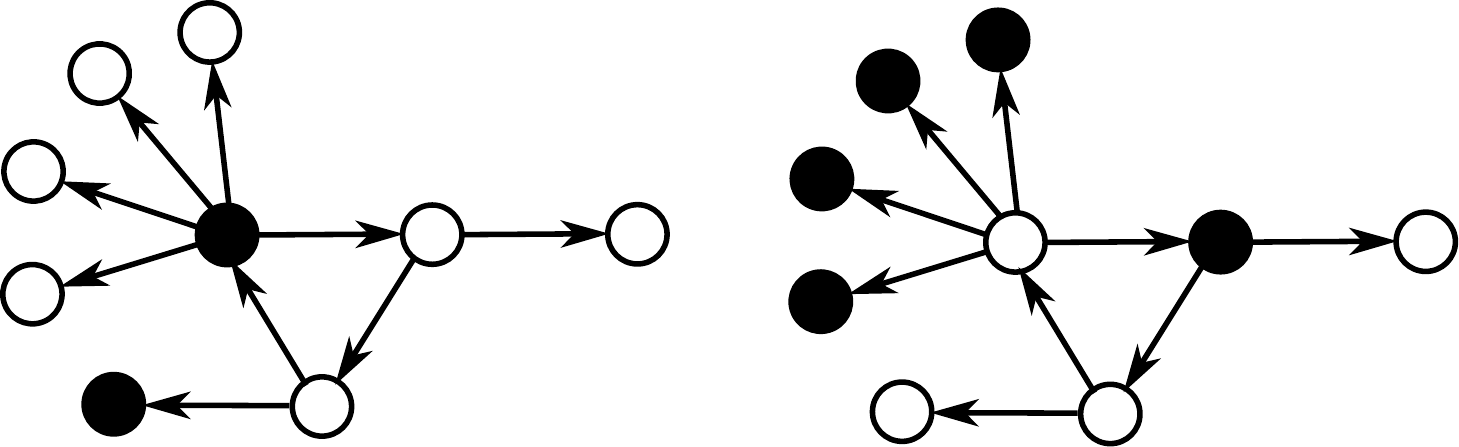}
\end{center}
\caption{Two different quasi-kernels for a directed graph. The first is smaller than
  the desired bound $|V(D)|/2$ while the second is greater than $|V(D)|/2$. Both 
  quasi-kernels can be constructed using Lemma \ref{lem:exists}.}
\label{fig:example}
\end{figure}

\section{Breakdown Sequences}
\label{sec:breakdown}
The purpose of this section is to investigate and formally define a sequence of 
breakdown steps for a directed graph. This will then be related to the 
construction of a specific type of quasi-kernel. For the proof, it will be necessary 
that such a construct covers the entire directed graph, as the premise of the graph 
being source-free is not inherited by strict subgraphs.  

The following lemma is adapted from \cite{chvatal}. It is based on a short induction 
proof for the existence of a quasi-kernel and is built around two types of steps
($\boldsymbol{\alpha}$ and $\boldsymbol{\beta}$) which are indicated in the proof for later reference.
\begin{lemma}
\label{lem:exists}
Every directed graph $D$ has a quasi-kernel.
\end{lemma}
\begin{proof}
By induction towards $|V(D)|$. If $V(D)\not=\emptyset$ then choose $v\in V(D)$
and thus, a quasi-kernel $Q'$ for $D-N^+[v]$ must exist by induction. If there exists
a $u\in Q'$ such that $u\longrightarrow v$ then $Q'$ is also a quasi-kernel
for $D$ ($\boldsymbol{\alpha}$-step). Otherwise, $Q'\cup\{u\}$ is a quasi-kernel for $D$ 
($\boldsymbol{\beta}$-step).
\end{proof}

The reader is encouraged to replicate the construction of both quasi-kernels
shown in Figure \ref{fig:example} by means of the induction proof in Lemma
\ref{lem:exists}. This will lead to the observation that the order in which 
vertices are chosen has direct implications for the size of the resulting 
quasi-kernel. Note that Lemma \ref{lem:exists} is sound, in the sense that
it always results in a quasi-kernel, but not complete, for it is not 
possible to construct every quasi-kernel in this way. However, the latter
will not be an obstacle for the method proposed in this work.

\begin{figure}
\begin{center}
\includegraphics[scale=.75]{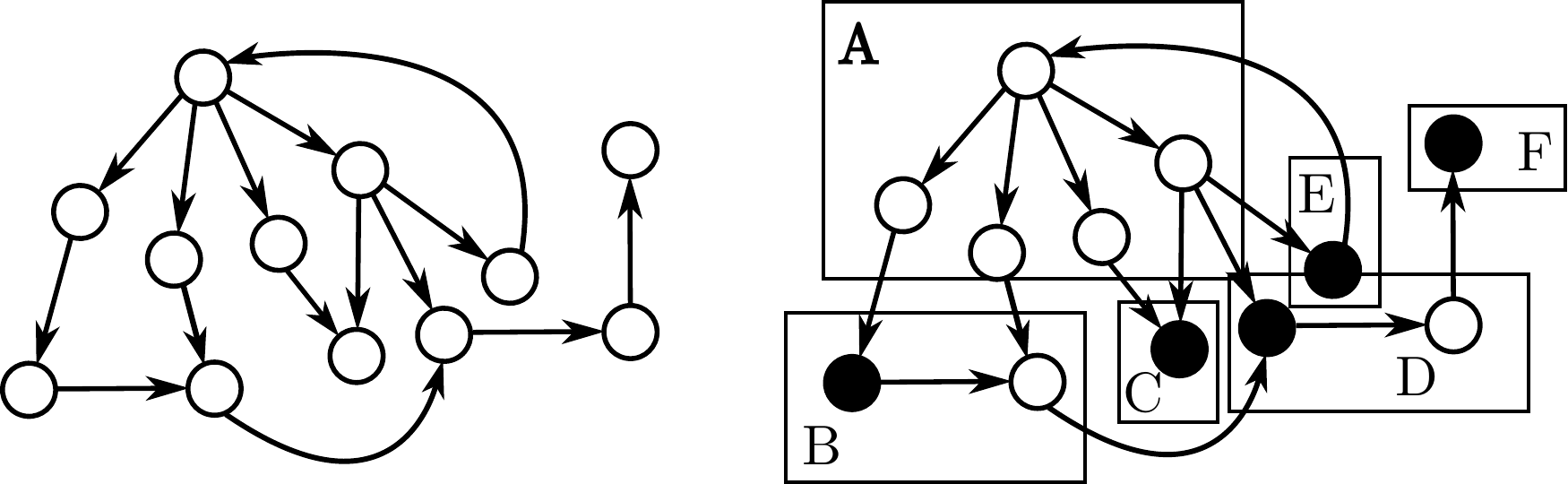}
\end{center}
\caption{A directed graph and a breakdown sequence in six steps A-F. A quasi-kernel
  as a subset of the first vertices from the breakdown sequence is shown in the
  second figure.}
\label{fig:breakdown}
\end{figure}

As shown in Figure \ref{fig:breakdown}, six recursive steps corresponding to 
Lemma \ref{lem:exists} are required to construct a quasi-kernel for this
particular directed graph. Step A is of type $\boldsymbol{\alpha}$, while 
steps B-F are of type $\boldsymbol{\beta}$. Note that the quasi-kernel in 
Figure \ref{fig:breakdown} is not minimal, as the vertex added in step F can 
be removed.

It is clear that if it is possible to apply Lemma \ref{lem:exists} in such a way 
that every $\boldsymbol{\beta}$-step is coupled to a non-empty neighborhood then this 
will result in a quasi-kernel of size at most $|V(D)|/2$. However, the directed 
graph in Figure \ref{fig:example} shows that this somewhat simplistic argument
will not work. Instead, an alternative counting method is applied.

In Section \ref{sec:main}, it will be shown that vertices in empty 
$\boldsymbol{\beta}$-steps must have emerged at some point as sources in
one of the recursive steps in the breakdown procedure in Lemma \ref{lem:exists}.
This leads to a distinction between $\boldsymbol{\beta}$-steps `under' an
$\boldsymbol{\alpha}$-step (e.g. step C in Figure \ref{fig:breakdown}),
and $\boldsymbol{\beta}$-steps `under' another $\boldsymbol{\beta}$-step
(e.g. step F in Figure \ref{fig:breakdown}. More formally, this notion of
`under' is defined as follows:

\begin{definition}
\label{def:under}
If $B$ is a breakdown sequence for a directed graph $D$ and if $(w,\emptyset)$
is in $B$ then $(w,\emptyset)$ is \underline{under} another step $(u,N)$ in $B$ if 
there exists a $v\in N$ such that $v\longrightarrow w$.
\end{definition}

Described abstractly, the counting method proposed in this paper works as follows:

\begin{enumerate}
\item Due to the fact that all non-empty $\boldsymbol{\beta}$-steps must have 
      disjoint out-neighborhoods, there must be at least twice as many vertices
      in these $\boldsymbol{\beta}$-steps compared to the number of 
      root-vertices.
\item It will be shown that $\boldsymbol{\beta}$-steps having an empty out-neighborhood, 
      which are itself under a non-empty $\boldsymbol{\beta}$-step, have an 
      \textsc{epon} in the set of root vertices of the $\boldsymbol{\alpha}$-steps.
\item Empty $\boldsymbol{\beta}$-steps which are under an $\boldsymbol{\alpha}$-step
      will be counted by means of the out-degree of the $\boldsymbol{\alpha}$-step
      they are under.
\end{enumerate}

Whether part 3 of the aforementioned enumeration is a valid part of the counting
method is open at present. Therefore, the following is postulated as an axiom here:

\begin{axiom}
\label{ax:axiom}
Every non-empty digraph $D$ has a vertex $u\in V(D)$ such that $d^+(u)>0$ and
$D-N^+[u]$ has at most $d^+(u)$ sources not present in $D$. 
\end{axiom}

A breakdown sequence is henceforth formally defined as:

\begin{definition}
\label{def:bds}
A \underline{breakdown sequence} (\textsc{bds}) for a directed graph $D$ is defined as a 
sequence $(u_1,N_1),(u_2,N_2),\ldots,(u_k,N_k)$, such that $u_i\in V(D)$ and $N_i\subseteq V(D)$ 
for all $1\leq i\leq k$, and is defined inductively by means of the following set 
of derivation rules:
\begin{enumerate}
\item If $V(D)=\emptyset$ then the empty sequence is a breakdown sequence for $D$.
\item If $A(D)=\emptyset$ and $u\in V(D)$ and if $B$ is a breakdown sequence for
      $D-\{u\}$ then $(u,\emptyset),B$ is a breakdown sequence for $D$.
\item If $u\in V(D)$ and $d^{+}(u)>0$ such that $D-N^+[u]$ introduces at most $d^{+}(u)$
      new sources then, if $B$ is a breakdown sequence for $D-N^+[u]$, then $(u,N^{+}(u)),B$
      is a breakdown sequence for $D$.
\end{enumerate}
\end{definition}

The correspondence between rule 3 in Definition \ref{def:bds} and Axiom \ref{ax:axiom}
should be clear to the reader. For additional clarity, if $B = (u_2,N_2),\ldots\,(u_k,N_k)$ is a 
breakdown sequence then the notation $(u_1,N_1),B$ is used to denote the 
concatenation of $(u_1,N_1)$ and $B$ (i.e. $(u_1,N_1),B = (u_1,N_1),(u_2,N_2),\ldots,(u_k,N_k)$).

The existence of a breakdown sequence is now immediate.

\begin{lemma}
\label{lem:bds}
Every directed graph $D$ has a breakdown sequence $B$.
\end{lemma}
\begin{proof}
Apply induction towards $|V(D)|$. If $A(D)\not=\emptyset$ then use Axiom \ref{ax:axiom}.
\end{proof}

For any breakdown sequence $B$, let $F(B)$ denote the set of vertices on the left-hand
side of the pairs in $B$. For example, if $B=(u_1,N_1),(u_2,N_2),\ldots,(u_k,N_k)$ is a 
\textsc{bds} then $F(B)=\{u_1,u_2,\ldots,u_k\}$. The class of quasi-kernels considered
here will be narrowed down by limiting these to quasi-kernels $Q\subseteq F(B)$ which
adhere to the following properties:

\begin{definition}
\label{def:qk_constr}
A quasi-kernel $Q$ is \underline{constructive} with regard to a breakdown sequence $B$ and a directed 
graph $D$ if for all $(w,N)$ in $B$ at least one of these conditions is met:
\begin{enumerate}
\item $w\in Q$, or
\item there exists a $v\in Q$ such that $v\longrightarrow w$, or
\item $N=\emptyset$ and there exist $u,v\in V(D)$ and $M\subseteq V(D)$ such that
      $(u,M)$ is in $B$ and $v\in M$ and $v\longrightarrow w$.
\end{enumerate}
\end{definition}

The first two conditions in Definition \ref{def:qk_constr} correspond to steps
of type $\boldsymbol{\beta}$ and $\boldsymbol{\alpha}$, respectively. In addition,
the third condition is required to enable minimality for $\boldsymbol{\beta}$-steps.
An example is the removal of vertices such as the one in step F in Figure \ref{fig:breakdown}.

Given a breakdown sequence $B$, it is possible to find a quasi-kernel which 
adheres to the properties listed in Definition \ref{def:qk_constr}.

\begin{lemma}
\label{lem:construct}
If $B$ is a \textsc{bds} for $D$ then there exists a constructive quasi-kernel
$Q$ for $B$ and $D$ such that $Q\subseteq F(D)$.
\end{lemma}
\begin{proof}
Apply induction towards the number of pairs in $B$, denoted here as $|B|$. Only
conditions 1 and 2 from Definition \ref{def:qk_constr} are used in the proof of
this lemma. Essentially, one can proceed in a very similar way to the setup of
Lemma \ref{lem:exists}. Note that the case for $|B|=0$ is handled immediately,
for $\emptyset$ is a constructive quasi-kernel if $V(D)=\emptyset$. Now, assume that 
$B=(u,N^+[u]),B'$ such that $B'$ is a breakdown sequence for $D-N^+[u]$. Further,
assume that $Q'\subseteq F(B')$ is a constructive quasi-kernel for $B'$ and 
$D-N^+[u]$ and distinguish between the following two cases:
\begin{enumerate}
\item If $(u,N^+[u]),B'$ results from rule 2 in Definition \ref{def:bds} then 
      $Q'\cup\{u\}$ is a constructive quasi-kernel for $B$ and $D$.
\item If $(u,N^+[u]),B'$ results from rule 3 in Definition \ref{def:bds} then 
      there are two possibilities. If there exists a $t\in Q'$ such that
      $t\longrightarrow u$ then $Q'$ satisfies condition 2 from Definition
      \ref{def:qk_constr} with regard to $u$. Hence, $Q'$ is a constructive
      quasi-kernel for $B$ and $D$. Otherwise, in case of absence of such a
      $t$, $Q'\cup\{u\}$ is a constructive quasi-kernel for $D$. 
\end{enumerate}
\end{proof}

\section{Main Theorem}
\label{sec:main}
This section completes the proof of the small quasi-kernel conjecture under
the assumption of Axiom \ref{ax:axiom}. This requires two additional
results. First, it will be shown that pairs $(u,\emptyset)$ in a breakdown
sequence at some point emerged as sources in the breakdown construction,
provided that the original graph was source-free (Lemma \ref{lem:rel_src}).
Subsequently, a result concerning \textsc{epon}s of vertices under a 
$\boldsymbol{\beta}$-step is required (Lemma \ref{lem:epon}).

\begin{lemma}
\label{lem:rel_src}
If $B$ is a \textsc{bds} for a directed graph $D$ and if $(u,\emptyset)$ is in 
$B$ such that $u$ is not a source in $D$ then $u$ emerged as a source in one of 
steps in the construction of $B$.
\end{lemma}
\begin{proof}
Apply induction towards the length of $B$, denoted here as $|B|$. The case for
$|B|=0$ is handled immediately as this clearly implies that $V(D)=\emptyset$.
Now, suppose that $B=(v,N^{+}[v]),B'$ for some $v\in V(D)$, such that $B'$
is a \textsc{bds} for $D-N^{+}[v]$. A distinction is made between the following
two cases:
\begin{itemize}
\item If $u\in V(D-N^{+}[v])$ then there are two possibilities. If $u$ is a
      source in $D-N^{+}[v]$ then $u$ emerged as a source in this step. 
      Otherwise, if $u$ is not a source in $D-N^{+}[v]$ then the result follows
      by induction.
\item If $u\in N^{+}[v]$ then $u=v$ and $N=\emptyset$ must hold due to the premise 
      that $(u,\emptyset)$ is in $B$. This means that $(v,N^{+}[v])$ and $B'$ were 
      not joined together by rule 3 in Definition \ref{def:bds}. This leaves rule 2
      as the only other possibility with the implication $A(D)=\emptyset$. This
      contradicts the premise that $u$ was not a source in $D$.
\end{itemize}
\end{proof}

\begin{lemma}
\label{lem:epon}
If $Q\subseteq F(B)$ is the smallest possible constructive quasi-kernel with regard to a
\textsc{bds} $B$ and a directed graph $D$, and if $(w,\emptyset)$ is in $B$ and $w\in Q$
and if $(w,\emptyset)$ is under a $\boldsymbol{\beta}$-step in $B$, then $w$ must have an
\textsc{epon} in $F(B)-Q$.
\end{lemma}
\begin{proof}
Suppose, towards a contradiction, that $w$ does not have an \textsc{epon} in $F(B)-Q$. It
will be shown that $Q-\{w\}$ is a constructive quasi-kernel with regard to $B$ and $D$.

Observe that, due to the premise that $Q$ is constructive, it is not possible for $w$ to
have an \textsc{epon} in $V(D)-F(B)$. Therefore, and due to the premise that $(w,\emptyset)$
is under a $\boldsymbol{\beta}$-step in $B$, $Q-\{w\}$ must be a quasi-kernel for $D$.

It remains to be shown that $Q-\{w\}$ is constructive with regard to $B$ and $D$. Clearly,
the removal of $w$ from $Q$ implies that $(w,\emptyset)$ is not longer covered by
condition 1 in Definition \ref{def:qk_constr}. Instead, $(w,\emptyset)$ is covered
by condition 3 in Definition \ref{def:qk_constr}.

If there exists a pair $(x,N)$ in $B$ such that $w\longrightarrow x$ and if $w$ is the 
only vertex in $Q$ which satisfies condition 2 in Definition \ref{def:qk_constr}, then
$x$ must be an \textsc{epon} in $F(B)-Q$ with regard to $w$, a contradiction.
\end{proof}

Given the previous results, it is now possible to prove the small quasi-kernel
conjecture under the assumption of Axiom \ref{ax:axiom}.

\begin{theorem}
Every source-free directed graph $D$ has a quasi-kernel $Q\subseteq V(D)$ such that $|Q|\leq |V(D)|/2$.
\end{theorem}
\begin{proof}
First, Lemma \ref{lem:bds} is used to derive the existence of a breakdown sequence $B$ for 
$D$ and then, by Lemma \ref{lem:construct}, there must exist a constructive quasi-kernel $Q\subseteq F(B)$. 
Let $Q$ be the smallest possible quasi-kernel which conforms to this property, 
for the given breakdown sequence $B$.

Due to the premise that $D$ is source-free and by Lemma \ref{lem:rel_src} it is known that 
elements $(u,\emptyset)$ in $B$ must have emerged as a source at some point in the
construction of $B$. Therefore, the following three pairwise disjoint sets can be defined:
\begin{center}
\begin{math}
\begin{array}{lcl}
R & = & \{(u,N)\mid u\in Q,\,N\not=\emptyset,\,(u,N)\,\,\textrm{in}\,\,B\} \\[8pt]
S & = & \{(u,\emptyset)\mid u\in Q,\,(u,\emptyset)\,\,\textrm{in}\,\, B,\,\,
	(u,\emptyset)\,\,\textrm{is under a}\,\,\boldsymbol{\beta}\textrm{-step in}\,\,B\} \\[8pt]
T & = & \{(u,\emptyset)\mid u\in Q,\,(u,\emptyset)\,\,\textrm{in}\,\, B,\,\,
	(u,\emptyset)\,\,\textrm{is under an}\,\,\boldsymbol{\alpha}\textrm{-step in}\,\,B\} - S\\
\end{array}
\end{math}
\end{center}
Clearly, it follows that $|Q|=|R\cup S\cup T|$. Due to the construction for
$B$ in Definition \ref{def:bds}, it is immediate that the neighborhoods 
in $B$ are disjoint. By Lemma \ref{lem:epon}, it is clear that the half-bound is 
satisfied for $S$ as all vertices have an \textsc{epon} in $F(B)-Q$. Furthermore, 
the number of elements in $T$ is dominated by the sum of the out-degrees of the
$\boldsymbol{\alpha}$-steps these are under, all of which correspond to disjoint
out-neighborhoods in $B$. The combination of these three disjoint cases leads to 
$|Q|\leq |V(D)-Q|$, which completes the proof.
\end{proof}

As a concluding remark, it is clear that if Axiom \ref{ax:axiom} can be shown 
to hold for every directed graph then this would complete the proof of the
small quasi-kernel conjecture. Therefore, it is postulated here as a conjecture.

\begin{conjecture}
\label{con:axiom}
Every non-empty digraph $D$ has a vertex $u\in V(D)$ such that $d^+(u)>0$ and
$D-N^+[u]$ has at most $d^+(u)$ sources not present in $D$. 
\end{conjecture}

Alternatively, another possible approach for the small quasi-kernel conjecture would 
be to combine the formulation proposed by Kostochka et al. \cite{kostochka} with a
method similar to the one proposed in this work. This will require further research.

\section*{Acknowledgments}
I would like to thank Peter Erd\H{o}s for pointing me to reference \cite{fete}
and I further thank Alexandr Kostochka for his words of encouragement to pursue
further investigation of this most interesting problem.

\end{document}